\documentclass[preprint,12pt]{elsarticle}
\usepackage{tikz, adjustbox}
\usepackage{xcolor}
\usepackage{xurl}
\usepackage[normalem]{ulem}

\newcommand{\ddt
}[1]{\frac{\partial {#1}}{\partial t}}

\usepackage{amssymb}
\usepackage{amsmath}
\usepackage{subcaption}

\journal{arXiv}

\begin{document}
\newpageafter{author}
\begin{frontmatter}

\title{Domain Decomposition-Based Coupling of High-Fidelity Finite Element and Reduced Order Operator Inference Models Using the Schwarz Alternating Method}

\author[VT,SandiaNM]{Ian Moore\corref{cor1}} \ead{ianm9123@vt.edu}
\cortext[cor1]{Corresponding Author}
\author[SandiaNM]{Anthony Gruber}
\author[SandiaNM]{Chris Wentland}
\author[SandiaNM]{Irina Tezaur}

\affiliation[SandiaNM]{organization = {Sandia National Laboratories$^1$},city = {Albuquerque},state = {NM},country={USA}}
\affiliation[VT]{organization={Virginia Tech Department of Mathematics},
city={Blacksburg}, state={VA},country={USA}}

\fntext[disclosure]{Sandia National Laboratories is a multimission laboratory managed and operated by National Technology \& Engineering Solutions of Sandia, LLC, a wholly owned subsidiary of Honeywell International Inc., for the U.S. Department of Energy’s National Nuclear Security Administration under contract DE-NA0003525. This work was funded through Sandia’s Laboratory Directed Research and Development program. The views expressed in the article do not necessarily represent the views of the U.S. Department of Energy or the United States Government.}

\begin{abstract}
 We propose a novel hybrid domain decomposition method that couples sub-domain-local high-fidelity finite element (FE) models with reduced order models (ROMs) using the Schwarz alternating method. By integrating the noninstrusive Operator Inference (OpInf) ROM, our approach accelerates the Schwarz process while allowing for geometry and mesh flexibility. We demonstrate the effectiveness of the new OpInf-FE method on a convection-dominated convection-diffusion-reaction problem, achieving stable and accurate predictive solutions while improving the ROM training process. 

\end{abstract}

\end{frontmatter}

\section{Introduction}
\label{intro}

Despite advancements in algorithms and high-performance computing (HPC), achieving detailed high-fidelity finite element (FE) analyses for complex multi-scale and multi-physics systems remains challenging, particularly for multi-query tasks like design optimization and uncertainty quantification. This is primarily due to two factors: (i) the lengthy process of creating high-quality meshes, identified as ``the single biggest bottleneck in [modeling and simulation] analyses"~\cite{SandiaLabNews}, and (ii) the long runtimes required for high-fidelity simulations of complex nonlinear processes. Although emerging data-driven reduced order models (ROMs) can alleviate runtime issues, they face challenges related to robustness, stability, accuracy, and refinement.

 We propose a novel method to address these challenges through hybrid domain decomposition (DD)-based models. Our approach couples subdomain-local high-fidelity models (HFMs) with subdomain-local ROMs using the Schwarz alternating method (SAM)~\cite{Schwarz:1870}. Building on previous work in solid mechanics~\cite{Mota:2017, Mota:2022} and computational fluid dynamics (CFD)~\cite{Wentland:2025}, we extend the SAM to couple HFMs with ROMs developed via the noninstrusive model reduction technique known as Operator Inference (OpInf)~\cite{PEHERSTORFERopinf}. 

The new OpInf-FE method simplifies mesh generation by enabling seamless coupling of non-conformal meshes and enhances the reliability of data-driven ROMs through online integration of HFM data. We extend our recent preliminary OpInf-Schwarz results~\cite{MooreCSRI} in a diffusion-dominated setting to a two-dimensional (2D) convection-dominated convection-diffusion-reaction (CDR) problem, achieving stable and accurate %
predictive solutions by combining low-order subdomain-local OpInf ROMs with a single subdomain-local HFM. Importantly, our method improves the OpInf ROM training process by removing the need for regularization in the least-squares minimization problem.
\section{The SAM and OpInf for a Model 2D CDR Problem}
We will now sketch the SAM coupling process for a model 2D CDR problem, given on the monolithic domain as, 
\begin{align}
    \ddt{u} - \epsilon \Delta u + \mathbf{b} \cdot \nabla u + \sigma u &= f, \,\,  &&t,\mathbf{x} \in (0, T] \times \Omega, \label{eq:CD_strong} 
\end{align}
with boundary data given as $g(t, \mathbf{x})$ on $\partial \Omega$ and an initial condition of 0 for simplicity. 
Above, $u : [0, T] \times \overline{\Omega} \rightarrow \mathbb{R}$ is the state, $\Omega \subset \mathbb{R}^2$ with closure $\overline{\Omega}$ is the spatial domain, $T \in \mathbb{R^+}$ is the final time, $\epsilon > 0$ and $\sigma \geq 0$ are %
parameters, $\mathbf{b} \in \mathbb{R}^2$ is the convection vector, 
and $f \in L^2(\Omega)$ is a forcing term. 

Problem~\eqref{eq:CD_strong} is challenging for ROMs to solve in the convection-dominated setting. While FE techniques can effectively solve~\eqref{eq:CD_strong}, these techniques require %
careful mesh generation, which can be one of the most expensive parts of HPC simulations on top of the already high cost of solving the FE model. The SAM has the potential to circumvent both of these issues through the construction of coupled hybrid DD-based OpInf-FE models.

To explain the basic SAM approach on equation~\eqref{eq:CD_strong}, consider a decomposition of $\Omega$ into two overlapping subdomains $\Omega_1$ and $\Omega_2$ (%
see Figure \ref{schematic:quadrants} for an idea of this), and express the linear operator governing \eqref{eq:CD_strong} as $\mathcal{L}: = (- \epsilon \Delta + \mathbf{b} \cdot \nabla + \sigma)$.  At simulation time $t_{n+1}>t_n$, the $(k+1)^{\mathrm{st}}$ Schwarz iteration is the sequential solution to the subdomain-local problems of~\eqref{eq:generic_schwarz_iter_Omega}, where $[S]$ is the indicator function for statement $S$ and $\Gamma_i$ are the $\Omega$-internal Schwarz
boundaries of $\Omega_i$ for $i = 1,2$. Equations \eqref{eq:generic_schwarz_iter_Omega}
advance the PDE in time with coupling only through Dirichlet boundary counditions, enabling flexible DD in space. Importantly, the SAM enables 
the use of %
nonconformal element types and solver methods 
in the subdomains $\Omega_1$ and $\Omega_2$.
\begin{equation} \label{eq:generic_schwarz_iter_Omega}
    \left \{
    \begin{array}{llll}
    &\ddt{}u_i^{(k+1)} + \mathcal{L} u_i^{(k+1)} = \,f, && \text{in } \Omega_i \times [t_n, t_{n+1}],\\
     &u_i^{(k+1)} = g, && \text{on } (\partial \Omega \cap \overline{\Omega}_i) \times [t_n,t_{n+1}], \\
    & u_i^{(k+1)} =  [i=1]\,u_2^{(k)} + [i=2]\,u_1^{(k+1)},  &&\text{on } \Gamma_i \times [t_n,t_{n+1}], 
    \end{array}
    \right.
\end{equation}

This work develops a new strategy for accelerating the Schwarz iteration process~\eqref{eq:generic_schwarz_iter_Omega} through assignment of both HFMs and noninstrusive OpInf ROMs to the subdomains comprising $\Omega$. This enables the use of ROMs in low-complexity regions alongside HFMs in more challenging regions, %
balancing the load according to the problem dynamics and minimizing the overall cost of computation.

OpInf \cite{PEHERSTORFERopinf} is a ROM technique that noninstrusively learns a set of reduced-order (and usually polynomial) equations mimicking the model form of the projection-based Galerkin ROM. This has the advantage of enabling order reduction without requiring intrusive access to the FOM code, but applying OpInf to coupled subdomains introduces new challenges. While a typical monolithic ROM often targets problems with constant-in-time boundary conditions (BCs) that can be eliminated through snapshot preprocessing, the non-homogeneous, $\Omega$-internal BCs of the SAM process \eqref{eq:generic_schwarz_iter_Omega}
vary spatially and temporally at every iteration. %

A FE approach to solving a problem like~\eqref{eq:CD_strong} with non-homogeneous boundaries is lifting, which solves a problem that can be written as: $\ddt{\mathbf{v}} = \mathbf{K} \mathbf{v} + \mathbf{B} \mathbf{g} + \mathbf{f}$. Here, $\mathbf{v} \in \mathbb{R}^{N}$ is a vector of FE coefficients %
for the interior of $\Omega$%
, $\mathbf{g} \in \mathbb{R}^m$ is the FE discretization of the boundary data on $\partial \Omega$, while $\mathbf{K} \in \mathbb{R}^{N \times N}$ and $\mathbf{B} \in \mathbb{R}^{N \times m}$ %
correspond to a FE discretization of~\eqref{eq:CD_strong}.  
To 
respect this structure with
OpInf, 
snapshots of $\mathbf{v}$ are collected at $n_t$ timesteps, a basis $\mathbf{\Psi} \in \mathbb{R}^{N \times r}$ is created using proper orthogonal decompositon, and compressed snapshot data $\mathbf{\Psi}^\intercal{\mathbf{v}_i}$ are obtained for $i =  1, \dots, n_t$. Once we have compressed data, we solve the low dimensional minimization problem, 
\begin{equation}
\label{eq:OpInf_2_norm}
    \underset{\widehat{\mathbf{K}}, \widehat{\mathbf{B}}, \widehat{\mathbf{f}}}{\text{arg\,min}}\sum_{i=1}^{n_t} \left\| \ddt{}\mathbf{\Psi}^\intercal{\mathbf{v}_i} - \widehat{\mathbf{K} }\mathbf{\Psi}^\intercal{\mathbf{v}_i} -  \widehat{\mathbf{B}} \mathbf{g}_i - \widehat{\mathbf{f}}\right\|_2^2 + \lambda^2 \left(\| \widehat{ \mathbf{K}} \|_F^2 + \| \widehat{\mathbf{B}} \|_F^2 + \| \widehat{\mathbf{f}} \|^2_2 \right). 
\end{equation}
The operators $\widehat{\mathbf{K}} \in \mathbb{R}^{r \times r}$, $\widehat{\mathbf{B}} \in \mathbb{R}^{r \times m}$ (which imposes the Schwarz BCs) and constant vector $\widehat{\mathbf{f}}\in \mathbb{R}^{r}$ correspond to their FE counterparts. %
The stability-promoting Tikhonov-type regularization strategy involving $\lambda$~\cite{shaneCombustion} is flexible and effective~\cite{FarcasOpInf}, but finding the optimal $\lambda$ is a critical and time-consuming process for complex problems.

\section{Numerical Results}

We now detail results for the CDR equation~\eqref{eq:CD_strong} on the unit square $\Omega = [0,1]^2$ with the choices $T = 5.0$, $\epsilon = 10^{-2}$, $\sigma = 10^{-3}$, and $\mathbf{b} = [\cos(\pi/3), \ \sin(\pi/3)]^T$,
which induce a sharp boundary layer in the upper right corner. The exact solution is dynamic from $t = 0$ to approximately $t = 1$, at which point the solution reaches a steady state. For consistency, ROMs in all subdomains
are trained on data between $t = 0$ and $t = 0.5$, and must predict beyond the training set until $t = 1.0$.  Moreover, these ROMs are required to remain stable until $T = 5.0$.

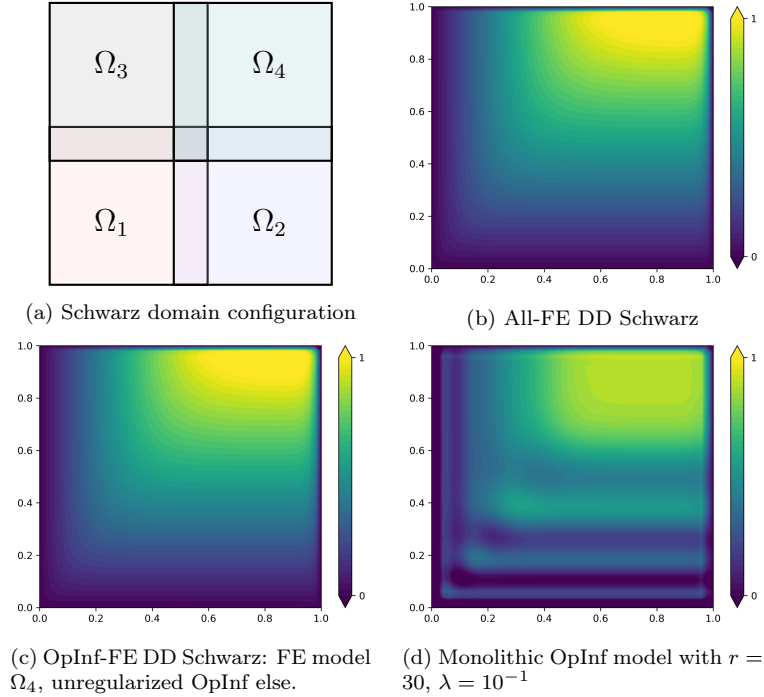
\begin{figure}[h]
    \centering
         \def\squareSize{2.1}
         \def\overlap{0.45}

     \begin{subfigure}[t!l]{0.35\textwidth}
         \centering
         \begin{tikzpicture}

         \def\squareplus{\squareSize + \overlap}
         \pgfmathsetmacro{\squareminus}{\squareSize - \overlap}

         \draw[fill=red!20, fill opacity=0.2,draw=black, thick] (0, 0) rectangle (\squareSize, \squareSize);
         \draw[fill=blue!20, fill opacity=0.2,draw=black, thick] ({\squareSize-\overlap}, 0) rectangle ({2*\squareSize-\overlap}, \squareSize);
         \draw[fill=darkgray!40, fill opacity=0.2,draw=black, thick] (0, {\squareSize-\overlap}) rectangle (\squareSize, {2*\squareSize-\overlap});
         \draw[fill=teal!40, fill opacity=0.2,draw=black, thick] ({\squareSize-\overlap}, {\squareSize-\overlap}) rectangle ({2*\squareSize-\overlap}, {2*\squareSize-\overlap});

         \node at ({\squareminus/2}, {\squareminus/2}) {$\Omega_1$};
         \node at ({\squareSize + \squareminus/2}, {\squareminus/2}) {$\Omega_2$};
         \node at ({\squareminus/2}, {\squareSize + \squareminus/2}) {$\Omega_3$};
         \node at ({\squareSize + \squareminus/2}, {\squareSize + \squareminus/2}) {$\Omega_4$};

         \end{tikzpicture}
         \caption{Schwarz domain configuration}
         \label{schematic:quadrants}
     \end{subfigure}
~
     \begin{subfigure}[t!r]{0.35\textwidth}
     \includegraphics[width = \textwidth]{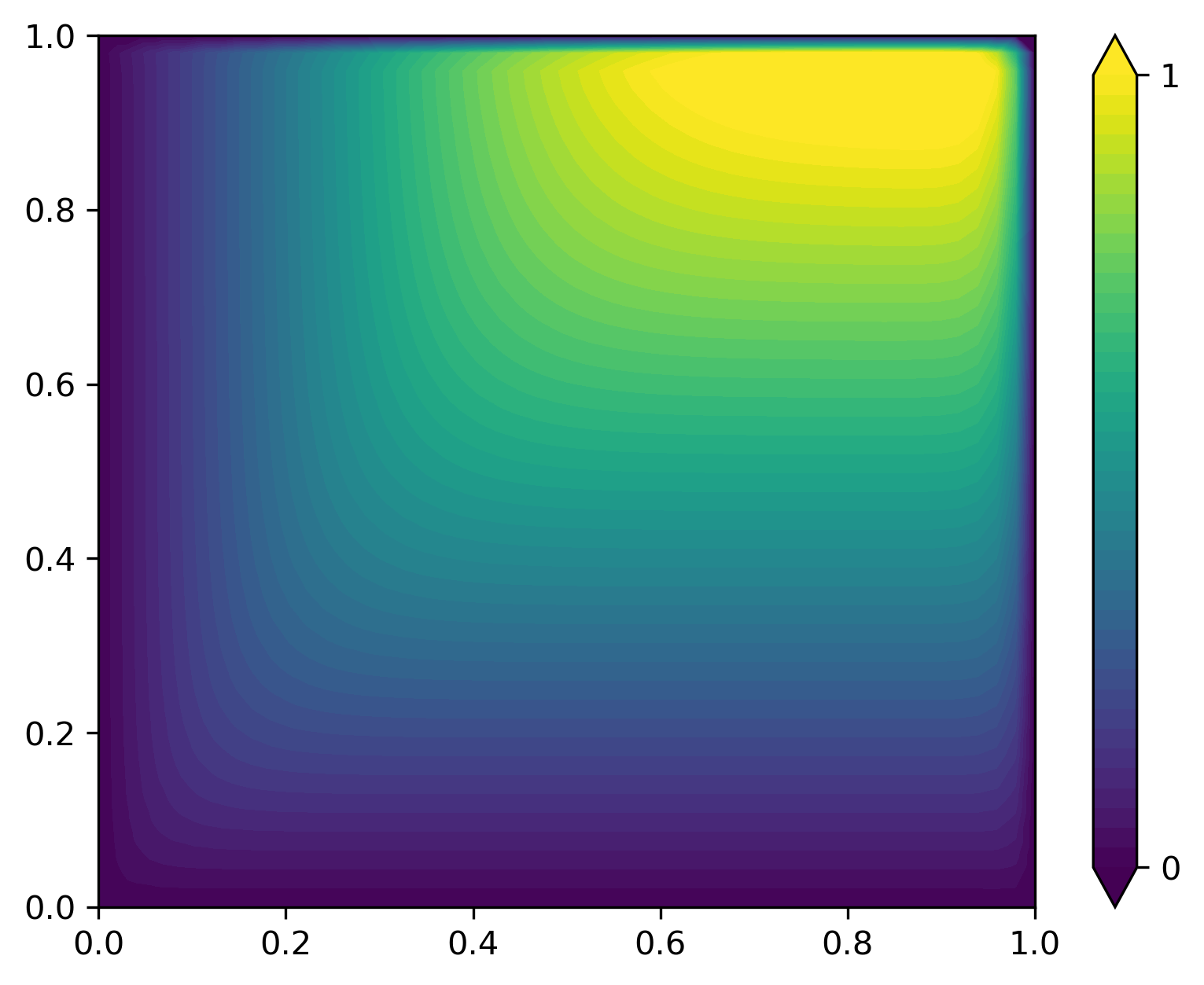}
     \caption{All-FE DD Schwarz } 
     \label{fig:All_FOM}
     \end{subfigure}\\
     
          \begin{subfigure}[t!]{0.35\textwidth}
         \includegraphics[width = \textwidth]{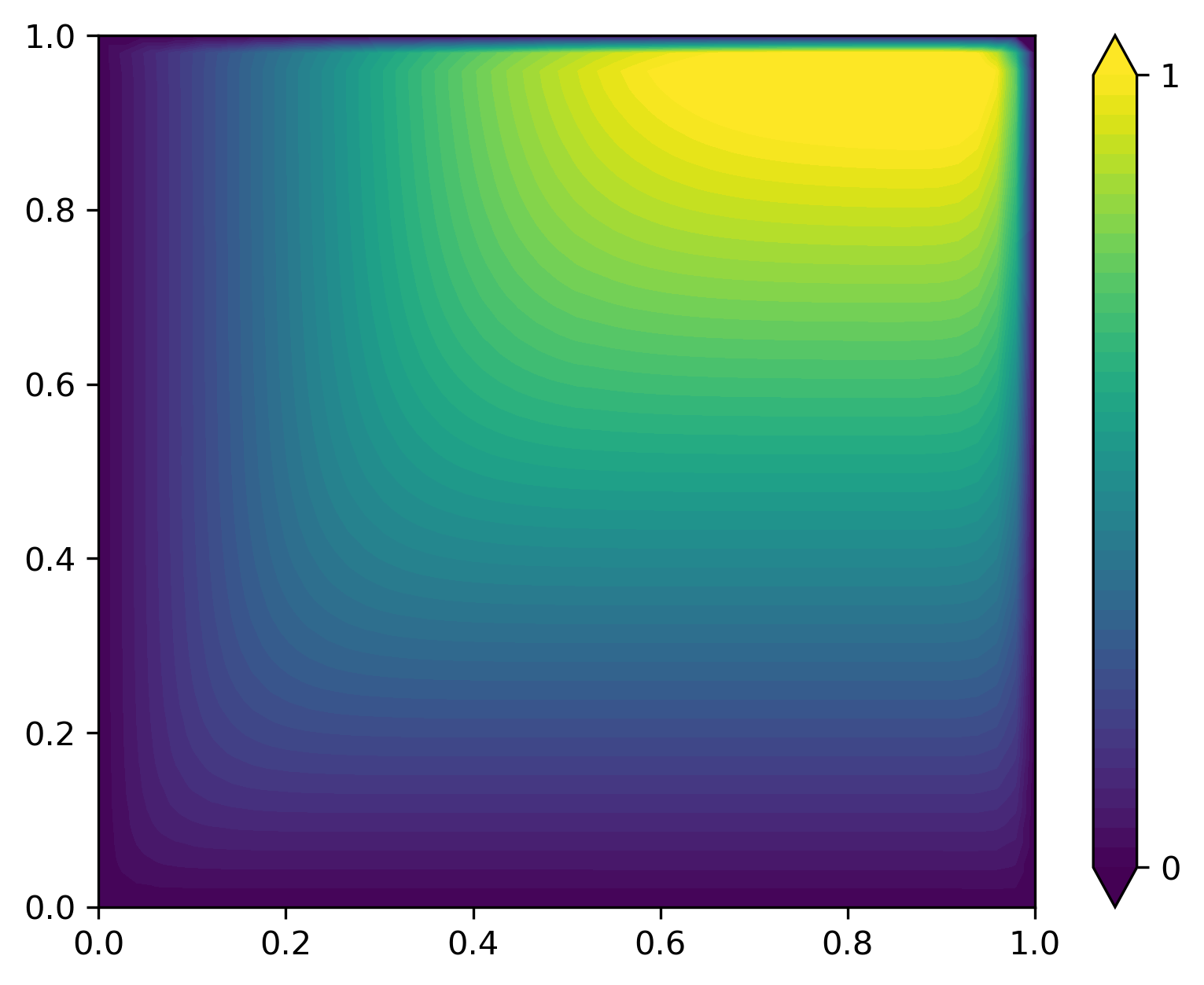}
         \caption{OpInf-FE DD Schwarz: FE model $\Omega_4$, unregularized OpInf else.}
         \label{fig:1_FOM}
     \end{subfigure}
     ~
     \begin{subfigure}[t!]{0.35\textwidth}
         \includegraphics[width = \textwidth]{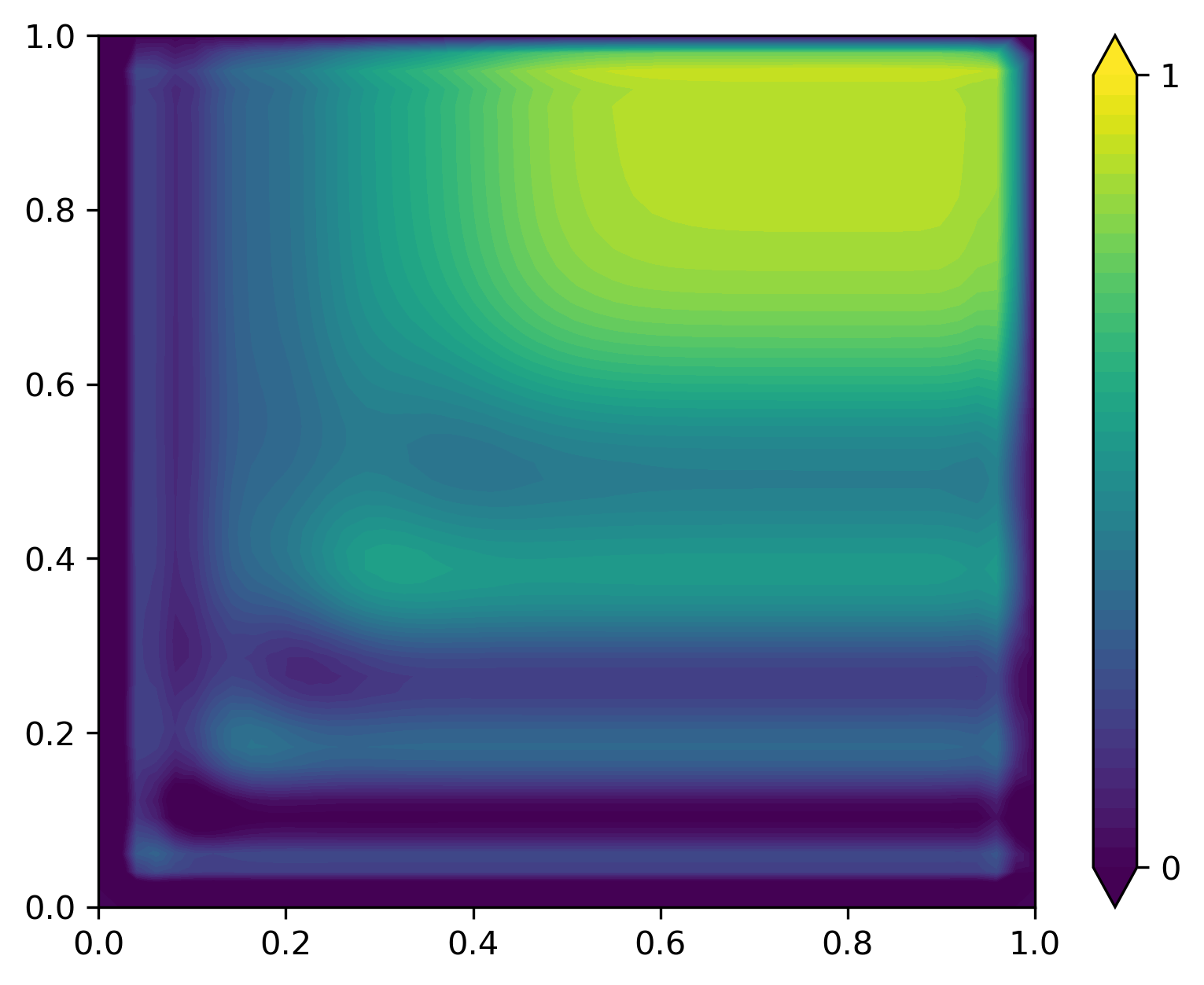}
         \caption{Monolithic OpInf model with $r = 30$, $\lambda = 10^{-1}$}
         \label{fig:mono_ROM}
     \end{subfigure}
        \caption{CDR numerical solutions at $T = 5.0$ %
        }
        \label{fig:comparison}
\end{figure}
A consistent DD of $\Omega$ (Figure~\ref{schematic:quadrants}) is used for all SAM simulations with underlying FE mesh resolution $h = 1/50$. First, the all-FE coupled Schwarz (Figure \ref{fig:All_FOM}) is the most accurate and expensive. Compare this to the OpInf-FE coupled Schwarz (Figure \ref{fig:1_FOM}), where 3 \textit{unregularized} OpInf models ($\lambda = 0)$ of basis dimension $r=10$ are coupled to a FE model in the boundary region $\Omega_4$, %
which produces a very similar result to the all-FE model without requiring additional regularization. Lastly, the monolithic OpInf model (Figure \ref{fig:mono_ROM}) struggles to produce an accurate solution even with a large basis size of $r = 30$ and strong regularization $(\lambda = 10^{-1})$ after tuning $\lambda$ in $[0, 1]$.%

\begin{table}[h]
    \centering
    \begin{tabular}{c|c|c|c}
      & All-FE Schwarz & OpInf-FE Schwarz & Mono. OpInf   \\\hline
       CPU time (s) &  18.7 & 8.9 & 0.2 \\
       Error & $\mathcal{O}(10^{-14})$ & $\mathcal{O}(10^{-4})$ & $\mathcal{O}(10^{-1})$
    \end{tabular}
    \caption{Simulation time and errors of various models (compared to monolithic FE model)}
    \label{tab:errors}
\end{table}

In Table~\ref{tab:errors}, we list the error and CPU time of the three models. ``Error" here refers the time-averaged relative pointwise error of each model in comparison to a monolithic FE model from $t = 0$ to $T = 5$. The Opinf-FE Schwarz model is twice as fast as the all-FE model, and is much more accurate than the monolithic Opinf model. %

\section{Conclusions \& Future Work}

We present a new 
DD-based method for coupling subdomain-local OpInf ROMs and FE simulations using the Schwarz alternating method.
While monolithic ROMs often struggle to resolve the varying spatial features inherent in convection-dominated problems, the proposed approach flexibly assigns OpInf and FE model types to regions where they individually excel. Through a CDR boundary value problem, we numerically demonstrate 
that our DD-based approach can remove the need for regularizing the OpInf minimization problem, thereby enabling the deployment of effective OpInf models with less tuning.  Beyond the superior accuracy of this hybrid OpInf-FE method over monolithic OpInf, we also demonstrate
significant CPU time savings compared to all-FE model couplings.

We emphasize that the proposed methodology is extensible to more challenging problems, for which the SAM displays additional benefits. In three-dimensional simulations, meshing complex geometries can itself be a daunting task, and the Schwarz framework allows for reuse of meshes as well as the alteration of mesh properties (such as resolution) in different subdomains. Future research directions 
include further easing meshing requirements through the development of the non-overlapping SAM and implementing parallel Schwarz to maximize time savings.%

\clearpage

\bibliographystyle{elsarticle-num} \bibliography{CDR_Opinf-Schwarz}

\end{document}